\documentclass[11pt]{amsart}
\usepackage[margin=2.6cm]{geometry}
\usepackage{amssymb,amsmath,booktabs,mathrsfs}
\usepackage[colorlinks=true,linkcolor=blue,citecolor=blue,urlcolor=blue]{hyperref}
\newtheorem{theorem}{Theorem}
\newtheorem{lemma}[theorem]{Lemma}
\newtheorem{corollary}[theorem]{Corollary}
\theoremstyle{remark}
\newtheorem{remark}[theorem]{Remark}

\newcommand{\Q}{\mathbb{Q}}
\newcommand{\Z}{\mathbb{Z}}
\newcommand{\F}{\mathbb{F}}
\newcommand{\sq}{\square}

\title[Torsion $\Z/5\Z$: exact rank four and rank at least five]{Elliptic curves with torsion $\Z/5\Z$: exact rank four over $\Q(k)$ and rank at least five infinitely often}
\author[A. Dujella]{Andrej Dujella}
\address{Department of Mathematics\\Faculty of Science\\University of Zagreb\\Bijeni{\v c}ka cesta 30, 10000 Zagreb, Croatia}
\email[A. Dujella]{duje@math.hr}
\author[J. C. Peral]{Juan Carlos Peral}
\address{Departamento de Matem\'aticas\\Universidad del Pa\'is Vasco\\Aptdo. 644, 48080 Bilbao, Spain}
\email[J. C. Peral]{juancarlos.peral@ehu.es}

\begin{document}

\begin{abstract}
We give two elliptic curves over $\Q(k)$ with torsion group $\Z/5\Z$ and rank exactly $4$, with four explicit
independent points. The upper bounds come from reducing the relevant elliptic surfaces modulo suitable primes and are
unconditional.
They match the known lower bound for the generic rank with this torsion, due to Eroshkin, and are not isomorphic to
Eroshkin's curve or to each other: their parametrising maps to $X_1(5)$ have degree $8$ instead of $12$, the
constructions are dihedral rather than biquadratic, and their branch loci differ. In each case two of the points come
from a ``doubling'' step in the spirit of Jaswal. Adding one further condition to the second curve, coming from a
quadratic twist over an intermediate conic, gives a fifth independent point over a genus-one curve of rank two. Hence there
are infinitely many elliptic curves over $\Q$ with torsion $\Z/5\Z$ and rank at least $5$, improving the known bound
$C(\Z/5\Z)\ge4$ to $C(\Z/5\Z)\ge5$.
\end{abstract}

\maketitle

\section{Results}

Every elliptic curve over $\Q$ with a point of order $5$ is a fibre of the universal curve
\[
E_t:\ y^2+(1-t)\,xy-t\,y=x^3-t\,x^2 ,
\]
with the $5$-torsion point $P=(0,0)$; its nonzero multiples have $x$-coordinates $0$ and $t$. As an elliptic surface over
the $t$-line it is rational, with fibres $I_5$ at $t=0,\infty$ and $I_1$ at the roots of $t^2-11t-1$, so $E(\Q(t))$
has rank $0$. The previously known lower bounds in Dujella's table \cite{Dtab} were $G(\Z/5\Z)\ge4$ for curves over $\Q(t)$
(Eroshkin, 2020) and $C(\Z/5\Z)\ge4$ for infinitely many curves over $\Q$ (Eroshkin, 2009).

For a point with $x$-coordinate $x$ on $E_t$, the $y$-coordinate is rational exactly when the discriminant
\[
\mathscr D_t(x)=\bigl((1-t)x-t\bigr)^2+4\bigl(x^3-tx^2\bigr)
\]
is a square.

\begin{theorem}[First family]\label{thm:A}
Let
\[
u=-\frac{20\,(k^2+1)^2}{(2k^2+5k-2)(5k^2-8k-5)},\qquad t=\frac{u^2-4}{3(u^2-1)} .
\]
Then $E_t$ over $\Q(k)$ has torsion group $\Z/5\Z$ and rank exactly $4$; points whose $x$-coordinates are
\[
\begin{aligned}
x_1&=-\frac{t}{3t-1}, & x_2&=-\frac{5153632t^3-3843928t^2+1025729t-1024}{1025^2\,(3t-1)},\\[1ex]
x_3&=\frac{u^2-4}{9(u-1)^2}, & x_4&=\frac{u^2-4}{9(u+1)^2}
\end{aligned}
\]
generate a subgroup of finite index in $E_t(\Q(k))$.
\end{theorem}

\begin{theorem}[Second family]\label{thm:B}
Let
\[
u=\frac{3k^4+16k^3-18k^2-16k+3}{3(k^2+1)^2},\qquad t=\frac{5-9u^2}{90(u^2-1)} .
\]
Then $E_t$ over $\Q(k)$ has torsion group $\Z/5\Z$ and rank exactly $4$; points whose $x$-coordinates are
\[
\begin{aligned}
x_1&=-\frac{8t+1}{20(10t+1)}, & x_2&=\frac{t(27t+4)}{4(8t+1)},\\[1ex]
x_3&=\frac{P(u)}{21600\,(u-1)(3u-5)^2(3u+5)(3u^2+4u+5)^2}, & x_4&=x_3\big|_{u\mapsto -u}
\end{aligned}
\]
generate a subgroup of finite index in $E_t(\Q(k))$, where
\[
\begin{aligned}
P(u)={}&19683u^9-715149u^8-323676u^7+4481892u^6+6257250u^5\\
&-9618750u^4-7357500u^3+9562500u^2-10078125u+6171875 .
\end{aligned}
\]
\end{theorem}

\begin{theorem}[Rank five over an elliptic curve]\label{thm:C}
Let $C'$ be the curve
\[
C':\ w^2=-5\,(3k^4+116k^3-218k^2-116k+3),
\]
which has genus one and the rational point $(k,w)=(-23,1850)$; its Jacobian is $y^2=x^3+x^2-73x+183$, of conductor
$20800$ and rank $2$. Over the function field $\Q(C')=\Q(k,w)$, the curve $E_t$ with $t=t(k)$ as in
Theorem~\ref{thm:B} has torsion group $\Z/5\Z$ and rank at least $5$: besides the four points of Theorem~\ref{thm:B} it
has the point with
\[
x_5=-\frac{(s^2-3s-9)\,P_7(s)}{3(s-6)^2(s+3)^2(s+6)(2s+3)^2(s^2+2s-39)},\qquad s=\frac{6(k^2+k-1)}{k^2-4k-1},
\]
where $P_7(s)=2s^7+21s^6-54s^5-1161s^4-1701s^3+15552s^2+45198s+30618$.
\end{theorem}

\begin{corollary}\label{cor:C}
There are infinitely many elliptic curves over $\Q$, with infinitely many distinct $j$-invariants, whose torsion group is
$\Z/5\Z$ and whose rank is at least $5$. Thus $C(\Z/5\Z)\ge5$.
\end{corollary}

In Theorems~\ref{thm:A} and~\ref{thm:B}, $t$ is a rational function of $k$ of degree $8$. For instance, in Theorem~\ref{thm:B},
\[
t=-\frac{k^8+24k^7+32k^6-168k^5-50k^4+168k^3+32k^2-24k+1}{40k\,(6k^6+7k^5-42k^4-14k^3+42k^2+7k-6)} .
\]

\section{Lower bounds}\label{sec:lower}

The lower bounds use the following observation.

\begin{lemma}\label{lem:eig}
Let $L$ be a field, $K=L(r,s)$ a biquadratic extension with $r^2,s^2\in L$, and $E/L$ an elliptic curve. If
$Q_1,\dots,Q_4\in E(K)$ have infinite order and lie in the eigenspaces of $E(K)\otimes\Q$ for $\operatorname{Gal}(K/L)$
belonging to four distinct characters, then they are independent.
\end{lemma}

\begin{proof}
Eigenspaces for distinct characters intersect trivially and are linearly independent, and each $Q_i\otimes1$ is a
nonzero vector in its eigenspace.
\end{proof}

A point whose $y$-coordinate involves a square root $\sqrt{c}$ only through $2y+a_1x+a_3=\sqrt{c}\cdot w$ with
$w\in L$ is sent to its negative by the automorphisms changing the sign of $\sqrt c$; so it lies in the eigenspace of the
character attached to $c$.

\begin{proof}[Proof of rank $\ge4$ in Theorem~\ref{thm:A}]
Put $L=\Q(u)$, $h=u(41u-80)$, $\tilde h=u(41u+80)$, and $K=L(r,s)$ with $r^2=h$, $s^2=\tilde h$. The classes of $h$ and
$\tilde h$ in $L^\times/L^{\times2}$ are independent, so $K/L$ is biquadratic. With $D=(2k^2+5k-2)(5k^2-8k-5)$ and $u$ as
in the theorem,
\[
r=\frac{20(k^2+1)(9k^2+2k-9)}{D},\qquad s=\frac{20(k^2+1)(k^2-18k-1)}{D},
\]
so $K\subseteq\Q(k)$, and $K=\Q(k)$ since both have degree $4$ over $L$. Let $\delta=(3t-1)(3t-4)$ and
$q=-(1573t+36)(3t-1)$. Direct calculation gives
\[
\begin{aligned}
\mathscr D_t(x_1)&=\delta\Bigl(\frac{2t^2}{(3t-1)^2}\Bigr)^2, &
\mathscr D_t(x_2)&=q\Bigl(\frac{2F(t)}{1025^3(3t-1)^2}\Bigr)^2,\\
\mathscr D_t(x_3)&=h\Bigl(\frac{(u-2)^2(u+2)}{27(u-1)^3(u+1)}\Bigr)^2, &
\mathscr D_t(x_4)&=\tilde h\Bigl(\frac{(u-2)(u+2)^2}{27(u-1)(u+1)^3}\Bigr)^2,
\end{aligned}
\]
with $F(t)=271502704t^4-169959988t^3+15426651t^2-342512t-87296$, and
\[
\delta=\Bigl(\frac{3u}{u^2-1}\Bigr)^2,\qquad q=\Bigl(\frac{rs}{u(u^2-1)}\Bigr)^2 .
\]
Hence the four points lie in $E_t(K)$, in the eigenspaces of the trivial character, the product character, the
character of $r$ and the character of $s$ respectively.

At $k=1/2$ we have $u=125/31$ and $t_0=1309/4888$. The fibre has good reduction at $59$ and $71$, with
$\#E_{t_0}(\F_{59})=70$ and $\#E_{t_0}(\F_{71})=75$, so its torsion subgroup is exactly $\Z/5\Z$. The specialised
points have $x$-coordinates $1309/961$, $566247748463/1206156099122$, $1309/8836$ and $1309/24336$, none of which is $0$
or $t_0$; so they have infinite order, hence so do the four points over $\Q(k)$. Lemma~\ref{lem:eig} gives independence.
Finally, $E_t(\Q(k))_{\rm tors}$ injects into $E_{t_0}(\Q)_{\rm tors}\cong\Z/5\Z$ and contains $P$.
\end{proof}

\begin{proof}[Proof of rank $\ge4$ in Theorem~\ref{thm:B}]
Put $L=\Q(u)$, $h=2(3u+5)$, $\tilde h=2(5-3u)$, $K=L(r,s)$ with $r^2=h$, $s^2=\tilde h$; again $K/L$ is biquadratic. For
$u$ as in the theorem,
\[
r=-\frac{4(k^2+k-1)}{k^2+1},\qquad s=\frac{2(k^2-4k-1)}{k^2+1},
\]
so $K=\Q(k)$. Let $\delta=5(10t+1)(18t+1)$ and $p=(8t+1)(18t+1)$. Direct calculation shows that $\mathscr D_t(x_1)$,
$\mathscr D_t(x_2)$, $\mathscr D_t(x_3)$ and $\mathscr D_t(x_4)$ are respectively $\delta$, $p$, $h$ and $\tilde h$ times
squares in $\Q(t)$ or $L$, and that
\[
\delta=\Bigl(\frac{4u}{3(u^2-1)}\Bigr)^2,\qquad p=h\tilde h\Bigl(\frac{u}{15(u^2-1)}\Bigr)^2 .
\]
So the four points lie in $E_t(K)$, in the eigenspaces of the trivial character, the product character, the character
of $r$ and the character of $s$.

At $k=3/2$ we have $t_0=-31919/604800$. The fibre has good reduction at $11$ and $13$, with
$\#E_{t_0}(\F_{11})=15$ and $\#E_{t_0}(\F_{13})=20$, so its torsion subgroup is exactly $\Z/5\Z$, and none of the four
specialised $x$-coordinates is $0$ or $t_0$. The proof concludes as before.
\end{proof}

\section{Upper bounds}\label{sec:upper}

\subsection*{Decomposition of the rank}
In both theorems $\Q(u)=\Q(t)(\sqrt\delta)$ and $\Q(k)=\Q(u)(\sqrt h,\sqrt{\tilde h})$ is biquadratic over $\Q(u)$,
where $\tilde h=h(-u)$. Write $E^d$ for the quadratic twist $d\,w^2=\mathscr F_t(x)$ of the universal curve, where
$\mathscr F_t(x)=4x^3+(t^2-6t+1)x^2+2(t^2-t)x+t^2$, and $E_\varphi$ for the pullback of $E$ to the $u$-line. Let
\[
r_d=\operatorname{rank}E^d(\Q(t)),\qquad r_h=\operatorname{rank}E_\varphi^{\,h}(\Q(u)) .
\]
The universal curve is an extremal rational elliptic surface (trivial lattice of rank $2+4+4=10$), so
$\operatorname{rank}E(\Q(t))=0$ by the Shioda--Tate formula. Decomposing under $\operatorname{Gal}(\Q(k)/\Q(u))$ and then
under $\operatorname{Gal}(\Q(u)/\Q(t))$, and using that $h\tilde h$ is the pullback of the class $q$ (Theorem~\ref{thm:A})
or $p$ (Theorem~\ref{thm:B}), and that $u\mapsto-u$ fixes $t$ and interchanges $h$ and $\tilde h$, we get
\[
\operatorname{rank}E_t(\Q(k))=r_\delta+r_{c}+r_{c\delta}+2r_h ,
\]
where $c=q$ for Theorem~\ref{thm:A} and $c=5(8t+1)(10t+1)$ for Theorem~\ref{thm:B} (so that $c\delta\equiv p$ there). By
Section~\ref{sec:lower}, $r_\delta,r_{c\delta}$ (Theorem~\ref{thm:B}) or $r_\delta,r_c$ (Theorem~\ref{thm:A}), and $r_h$ are at least $1$.

\subsection*{Reduction of the surfaces}
For each twist, let $X$ be the corresponding elliptic surface over $\mathbb P^1$ (a K3 surface for the twists over the
$t$-line, with fibres $2I_5+2I_1+2I_0^*$; a surface with $\chi=3$ for the $h$-twists, with fibres $4I_5+4I_1+2I_0^*$). If
$\ell$ is a prime not dividing $2\cdot3\cdot5$ at which the branch points of the twist stay distinct from each other and
from the singular fibres, then $X$ has good reduction at $\ell$ with the same fibre configuration. Let
\[
L_\ell(T)=\det\bigl(1-T\,\mathrm{Frob}_\ell\mid H^1(\mathbb P^1_{\overline{\F}_\ell},j_*R^1f_*\Q_{\ell'})\bigr),
\]
of degree $4$ for the twists over the $t$-line and $8$ for the $h$-twists. Specialisation of N\'eron--Severi groups is
injective \cite{MP}, and so is the cycle class map on $\mathrm{NS}\otimes\Q_{\ell'}$. By the Shioda--Tate decomposition
\cite{SS}, the Mordell--Weil group tensored with $\Q_{\ell'}$ is the orthogonal complement of the trivial lattice
(spanned by the zero section, a fibre and the fibre components), and in $H^2(X_{\overline{\F}_\ell},\Q_{\ell'}(1))$ this
complement lies in the parabolic cohomology above, Tate-twisted. The trivial lattice is stable under Frobenius, so the
orthogonal projection onto its complement commutes with Frobenius. A section defined over the constant field $\Q$ has a
Galois-invariant class, whose specialisation is fixed by $\mathrm{Frob}_\ell$ on $H^2(1)$; so the Mordell--Weil group over
$\Q$ embeds, after $\otimes\Q_{\ell'}$, into the parabolic cohomology as classes on which Frobenius acts by $\ell$
(compare \cite{vL}). Hence
\[
r\le\operatorname{ord}_{T=1/\ell}L_\ell(T),
\]
unconditionally. By the Grothendieck--Lefschetz formula, $L_\ell(T)=\exp\bigl(\sum_{n\ge1}S_nT^n/n\bigr)$, where $S_n$ is
the sum of the local traces over $\mathbb P^1(\F_{\ell^n})$ ($\ell^n+1-\#$ points for good fibres, $\pm1$ for split or
non-split multiplicative fibres, $0$ at $I_0^*$). For the twist by $d$, the local trace at a fibre $t$ is
$\chi(d(t))\,a_t$, where $\chi$ is the quadratic character of $\F_{\ell^n}$ and $a_t$ is the local trace of $E_t$ (both
fibres $t=0,\infty$ of $E$ are split of type $I_5$); so one point count of $E_t$ over each field serves all twists. The polynomials are fixed by $S_1,\dots,S_4$ and the functional
equation $c_{d-i}=\varepsilon\,\ell^{d-2i}c_i$ ($d=\deg L_\ell$), up to the sign $\varepsilon$; when both signs are
possible, the sign giving multiplicity $0$ for an $h$-twist is excluded by the known section.

\subsection*{Data}
For Theorem~\ref{thm:A} all bounds come from $\ell=29$, for Theorem~\ref{thm:B} from $\ell=29,31,47$. In each case the
leading coefficient and discriminant of the twisting polynomial, its resultant with $t(t^2-11t-1)$ (respectively with the
polynomial $W(u)$ cutting out the singular fibres of $E_\varphi$), and $125=\operatorname{disc}(t^2-11t-1)$ are units at
$\ell$, and the branch points of $u\mapsto t$ lie over smooth fibres modulo $\ell$.

\medskip
{\small
\noindent\begin{tabular}{@{}lllll@{}}
\toprule
& component & $\ell$ & $S_1,S_2,S_3,S_4$ & $L_\ell(T)$; $\operatorname{ord}_{1/\ell}$\\
\midrule
Thm~\ref{thm:A} & $\delta$ & $29$ & $-43,-1849,28982,-27889$ & $(1-841T^2)(1-43T+841T^2)$; $1$\\
& $q$ & $29$ & $-43,-1849,28982,-27889$ & $(1-841T^2)(1-43T+841T^2)$; $1$\\
& $q\delta$ & $29$ & $0,-2204,0,400316$ & $1-1102T^2+707281T^4$; $0$\\
& $h$ & $29$ & $-13,27,-53098,-930053$ & $(1-841T^2)R_6(T)$; $1$\\
\midrule
Thm~\ref{thm:B} & $\delta$ & $29$ & $-13,-169,30602,-2289169$ & $(1-841T^2)(1-13T+841T^2)$; $1$\\
& $c$ & $47$ & $20,-5476,292820$ & $(1+47T)^2(1-74T+2209T^2)$; $0$\\
& $c\delta\equiv p$ & $31$ & $38,-1444,-54682$ & $(1-961T^2)(1+38T+961T^2)$; $1$\\
& $h$ & $31$ & $-80,-1954,-34280,113342$ & $(1-961T^2)R'_6(T)$; $1$\\
\bottomrule
\end{tabular}
}

\medskip
\noindent Here
\[
\begin{aligned}
R_6(T)&=1-13T+939T^2-29174T^3+789699T^4-9194653T^5+594823321T^6,\\
R'_6(T)&=1-80T+3184T^2-95480T^3+3059824T^4-73881680T^5+887503681T^6,
\end{aligned}
\]
with $R_6(1/29)=1800/841$ and $R'_6(1/31)=250/961$. For the $h$-twist of Theorem~\ref{thm:A} at $29$ the coefficients
$c_0,\dots,c_4$ are $1,-13,98,-18241,0$: the sign $\varepsilon=+1$ gives multiplicity $0$ (excluded) and $\varepsilon=-1$
gives the displayed polynomial. For Theorem~\ref{thm:B} at $31$ they are $1,-80,2223,-18600,0$, and again only the sign
with multiplicity $1$ is compatible with the known section. For the twists over the $t$-line, $S_3$ fixes the sign. The rows for $\delta$ and $q$ in Theorem~\ref{thm:A}
coincide; this is not a misprint.

\begin{proof}[Proof of rank $\le4$ in Theorems~\ref{thm:A} and~\ref{thm:B}]
The table gives $r_\delta+r_c+r_{c\delta}+2r_h\le1+1+0+2=4$ for Theorem~\ref{thm:A} and $\le1+0+1+2=4$ for
Theorem~\ref{thm:B}. Together with Section~\ref{sec:lower} this proves rank exactly $4$, and that the four points generate a subgroup of
finite index.
\end{proof}

\begin{remark}
The choice of primes matters. At $\ell=29$ the twists $c$ and $c\delta$ of Theorem~\ref{thm:B} have multiplicities $1$ and
$2$, and its $h$-twist has multiplicity $2$, so $\ell=29$ gives weaker upper bounds for these components. (The extra
eigenvalues may come from sections defined over a quadratic extension of the constant field, but we have not identified
such sections.) The primes $\ell=47$ and $\ell=31$ separate these contributions and give sharp bounds.
\end{remark}

\begin{remark}[Specialisations]\label{rem:spec}
If $k$ is $71$-integral with $k\equiv1/2\pmod{71}$, the curve of Theorem~\ref{thm:A} reduces to the curve with $t=23$
over $\F_{71}$, of order $75$; if $k$ is $11$-integral with $k\equiv3/2\pmod{11}$, the curve of Theorem~\ref{thm:B}
reduces to the curve with $t=4$ over $\F_{11}$, of order $15$. These orders are odd, so such fibres have no rational
$2$-torsion, and by Mazur's theorem their torsion is exactly $\Z/5\Z$. By Silverman's specialisation theorem \cite{Sil}
all but finitely many of them have rank at least $4$, and their $j$-invariants take infinitely many values. Several fibres
of the second family have rank $5$, for example $k=4,5,1/4,5/2$ (PARI's \texttt{ellrank} bounds; for curves without
rational $2$-torsion the upper bound may rely on GRH).
\end{remark}

\section{Rank five: proof of Theorem~\ref{thm:C}}\label{sec:five}

\subsection*{The fifth point}
In the setting of Theorem~\ref{thm:B}, $h\tilde h=4(25-9u^2)$, and the conic $25-9u^2=\sq$ is parametrised by
$u=5(s^2-9)/(3(s^2+9))$; on it $h\equiv5(s^2+9)$ modulo squares. Over the $s$-line the curve $E_t$ has rank $2$ (the
points $x_1,x_2$), and its twist by $h$ has rank $2$ (the points $x_3,x_4$). The new input is the twist by
\[
D(s)=-15(s^2+2s-39),
\]
a class on the conic that is not pulled back from the $u$-line; it carries the point $x_5$:
$\mathscr D_t(x_5)$ is $D(s)$ times a square in $\Q(s)$. With $s=6(k^2+k-1)/(k^2-4k-1)$, which lies in $\Q(k)$ because
$\sqrt{h\tilde h}\in\Q(k)$, one has $u(s)=u(k)$ and
\[
D(s)=\frac{9\,(-5)(3k^4+116k^3-218k^2-116k+3)}{(k^2-4k-1)^2},
\]
so $\sqrt{D}=3w/(k^2-4k-1)\in\Q(C')$.

\begin{proof}[Proof of Theorem~\ref{thm:C}]
All identities above are exact (checked in PARI/GP). Since $C'$ has genus one, $w\notin\Q(k)$ and $\Q(C')=\Q(k)(w)$ is
quadratic over $\Q(k)$. Hence $E_t(\Q(C'))\otimes\Q$ splits into the $(+1)$- and $(-1)$-eigenspaces of $w\mapsto-w$, which
are $E_t(\Q(k))\otimes\Q$ and the part coming from the twist by $w^2$. The points $x_1,\dots,x_4$ lie in the first and are
independent by Theorem~\ref{thm:B}; $x_5$ lies in the second, since $w\mapsto-w$ changes the sign of $\sqrt D$ and so
sends the point to its negative. The torsion of $E_t(\Q(C'))$ injects into that of a fibre: at $k=-23$, $s=303/62$ and
$t_0=239761/76507200$, with good reduction at $29$ and $31$ and $\#E_{t_0}(\F_{29})=35$,
$\#E_{t_0}(\F_{31})=30$, so the torsion is $\Z/5\Z$, generated by $(0,0)$ and fixed by $w\mapsto-w$. If $x_5$ were torsion
it would equal its own negative, hence be the identity, which it is not. So $x_5$ has infinite order, and the rank is at
least $5$.
\end{proof}

\begin{proof}[Proof of Corollary~\ref{cor:C}]
$C'$ has a rational point and its Jacobian has rank $2$, so $C'(\Q)$ is infinite. By Silverman's specialisation
theorem \cite{Sil}, all but finitely many of its points give curves of rank at least $5$. The quartic defining $C'$ has
unit leading coefficient and discriminant at $29$, and $29\nmid20800$, so $C'$ has good reduction at $29$. Take
$P_0=(-23,1850)$ as the origin of $C'$. Reduction modulo $29$ is then a homomorphism $C'(\Q)\to C'(\F_{29})$, whose kernel
has finite index and is therefore infinite. Every point in the kernel reduces to $\overline{P_0}=(6,23)$, so it has $k$
integral at $29$ and $k\equiv-23\pmod{29}$. For these points $t(k)\equiv2\pmod{29}$, and the fibre reduces to $E_2/\F_{29}$, which has $35$ points; so it has no
rational $2$-torsion, and by Mazur's theorem its torsion is exactly $\Z/5\Z$. Since $t$ is non-constant on $C'$, infinitely
many $j$-invariants occur.
\end{proof}

\subsection*{Fibres}
At several small points of $C'$ the five points are visibly independent, and some fibres have rank $6$ (PARI's
\texttt{ellrank} bounds; possibly GRH-dependent):

\medskip
\begin{center}
\begin{tabular}{llll}
\toprule
$k$ & $s$ & $\det$ of the $5\times5$ height matrix & rank bounds\\
\midrule
$-23$ & $303/62$ & $9.81\times10^{7}$ & $6,6$\\
$19/17$ & $-237/122$ & $1.35\times10^{8}$ & $6,6$\\
$-19/8$ & $174/181$ & $2.09\times10^{8}$ & $5,5$\\
\bottomrule
\end{tabular}
\end{center}

\subsection*{The exact rank}
Since $\Q(C')=\Q(k)(w)$,
\[
\operatorname{rank}E_t(\Q(C'))=4+r_D+r_{Dh},
\]
where $r_D$ and $r_{Dh}$ are the ranks over $\Q(s)$ of the twists by $D$ and $Dh$ of $E_t$ over the $s$-line; the $4$ is
exact by Theorem~\ref{thm:B}, $r_D\ge1$ and $r_{Dh}\ge0$. The corresponding elliptic surfaces have parabolic cohomology of
dimension $16$ and $20$, so the method of Section~\ref{sec:upper} needs point counts over $\F_{\ell^8}$ and $\F_{\ell^{10}}$.
The primes $7$ and $11$ are of bad reduction for these twists, and already $\ell=13$ is beyond our computation, so we have
not determined the exact rank. Nagao-type averages give $0.93$ for $D$ and $-0.19$ for $Dh$, and the fibre at $k=-19/8$
has rank $5$ by \texttt{ellrank} (possibly GRH-dependent); if the rank over $\Q(C')$ were at least $6$, this fibre would be one of the finitely many exceptional
fibres of Silverman's theorem. We therefore expect the rank to be exactly $5$.

\begin{remark}[Heuristics]
If the number of curves with torsion $\Z/5\Z$ and height at most $H$ is of order $H^{1/6}$ \cite{HS}, the heuristic of
Park, Poonen, Voight and Wood \cite{PPVW} suggests the rank bound $1+24/6=5$ for this torsion. Rank $5$ is then the borderline
case, in which the model allows $H^{o(1)}$ curves; the family of Corollary~\ref{cor:C}, parametrised by the rational points
of an elliptic curve, is very sparse ($O(\log H)$ curves of height at most $H$), so the result is consistent with
the heuristic.
\end{remark}

\section{How the curves arise}\label{sec:arise}

In both cases a positive conic $\delta(t)=\sq$ is parametrised by $u$, with deck involution $\iota(u)=-u$; let
$E_\varphi=E_{t(u)}$. A twist $h$ of $E_\varphi$ with one point gives a second point on the twist by
$\tilde h=h\circ\iota$, as in Jaswal's doubling lemma \cite{J}. In his setting $h$ and $\tilde h$ have four distinct
roots and the base has genus one. Here $h$ has a root at a fixed point of $\iota$ ($u=0$, respectively $u=\infty$), so
the curve $h=\sq$, $\tilde h=\sq$ is branched at only three points and is rational. The four points come from:

\medskip
\begin{center}
\begin{tabular}{lll}
\toprule
point & Theorem~\ref{thm:A} & Theorem~\ref{thm:B}\\
\midrule
$x_1$ & twist of $E$ by $\delta=(3t-1)(3t-4)$ & twist of $E$ by $\delta=5(10t+1)(18t+1)$\\
$x_2$ & twist of $E$ by $q=-(1573t+36)(3t-1)$ & twist of $E$ by $p=(8t+1)(18t+1)$\\
$x_3$ & twist of $E_\varphi$ by $h=u(41u-80)$ & twist of $E_\varphi$ by $h=2(3u+5)$\\
$x_4$ & twist of $E_\varphi$ by $\tilde h=u(41u+80)$ & twist of $E_\varphi$ by $\tilde h=2(5-3u)$\\
\bottomrule
\end{tabular}
\end{center}

\medskip
In both cases the class of $x_2$ equals $h\tilde h$ on the $u$-line, and it shares a root with $\delta$ on the $t$-line.
The positive twists were located by Nagao-type averages of fibre traces ($106$ primes between $300$ and $1000$) over all
$77$ positive conics $\delta$ with rational roots found in a twist search; up to the automorphism $t\mapsto-1/t$ of the
universal curve, these averages single out these two configurations as the only ones of total rank $4$. The sections were reconstructed by
interpolating $x$-coordinates from rank-$1$ specialisations and then verified symbolically.

\section{Comparison}

Eroshkin's curve over $\Q(T)$ \cite{Er} is $E_{-s(T)}$ in the notation above, with
\[
s(T)=\frac{(T^3+2T^2-4)(T^3-2T^2+4)(T^3-2T-2)(T^3-2T+2)}{12\,T^2(T^2-4)(T^2-1)(T^2+2T+2)(T^2-2T+2)} .
\]
We checked that his four points lie on the curve and are independent at $T=7/3$ (determinant $5922.44\ldots$). His map
$T\mapsto s(T)$ has degree $12$ and is invariant under $T\mapsto -T$ and $T\mapsto 2/T$, a biquadratic (Klein four)
structure. The maps $k\mapsto t$ of Theorems~\ref{thm:A} and~\ref{thm:B} have degree $8$, and their Galois closures over
the $t$-line are dihedral. The degree of the map to $X_1(5)$ is unchanged by reparametrisation and by $t\mapsto-1/t$,
and $12/8$ is not an integer, so neither of our curves is isomorphic to a base change of Eroshkin's or conversely.
An isomorphism between two curves of the form $E_t$ maps the rational point $P$ of order $5$ to one of $\pm P,\pm2P$;
on the parameter, $\pm P$ give the same value of $t$ and $\pm2P$ give $-1/t$. So after a reparametrisation of $k$ the two
maps $k\mapsto t$ would have to agree up to $t\mapsto-1/t$. The two new curves are not isomorphic to each other: the
branch values of the maps $k\mapsto t$ are
$\{1/3,4/3,-36/1573\}$ and $\{-1/10,-1/18,-1/8\}$, and these sets are not related by the identity or by $t\mapsto-1/t$.

\section{Searches}

The fifth point of Theorem~\ref{thm:C} was found by searching for quadratic twists of the rank-$2$ curves over the three
conics $h=\sq$, $\tilde h=\sq$ and $h\tilde h=\sq$ that lie under the rational bases of Theorems~\ref{thm:A}
and~\ref{thm:B}; only the conic $h\tilde h=\sq$ of Theorem~\ref{thm:B} produced a class not pulled back from the
$u$-line. The following searches detected nothing beyond rank $4$ over $\Q(k)$ (Nagao averages, so heuristic):
\begin{itemize}
\item quadratic twists of the curve of Theorem~\ref{thm:A} in its base parameter (coefficients up to $80$, and $a\le3$,
$|b|\le300$, $|c|\le1500$), linear twists (coefficients up to $250$) and products of linear factors (coefficients up to
$40$); a positive one with a soluble conic would give rank $5$ over $\Q(w)$;
\item the same searches on the rank-$2$ curves over the conics $h=\sq$ and $h\tilde h=\sq$ of Theorem~\ref{thm:A}
and $h=\sq$ of Theorem~\ref{thm:B}, apart from classes pulled back from the $u$-line;
\item twists $k(u-r)\,g(u)$ of $E_\varphi$ for which the base of Theorem~\ref{thm:A} or~\ref{thm:B} acquires a genus-one
double cover ($r$ of height up to $300$);
\item the same searches on Eroshkin's curve, where by its symmetry a positive twist $h$ would come with $h\circ g$ and
give rank $6$ over a genus-one curve;
\item a configuration as in Section~\ref{sec:arise} with a piece of rank $2$: none occurred among the $77$ conics.
\end{itemize}

\section{Computations and use of AI}

The searches, sections and specialisations were computed with PARI/GP 2.15.4 and Python/NumPy; all identities used in
the proofs, including those for Theorem~\ref{thm:C}, were checked exactly in PARI/GP. The trace sums $S_n$ in
Section~\ref{sec:upper} were computed by point counting over $\mathbb P^1(\F_{\ell^n})$ in PARI/GP; for Theorem~\ref{thm:A} they were
obtained independently in two separate implementations and agree. Scripts reproducing all identities, all trace sums
$S_n$ and the resulting $L$-polynomials, and the numerical data of Sections~\ref{sec:lower}--\ref{sec:five}, are
available as ancillary files. The exact-rank argument by reduction of the
surfaces follows the second author's note on Theorem~\ref{thm:A}. The construction was
found in an investigation carried out with Claude (Anthropic), an AI system, following ideas of Jaswal \cite{J}. The
eigenspace proof of Theorem~\ref{thm:A}, which needs no rank upper bounds, follows a suggestion received on an earlier
version of this note. The authors directed the work and take full responsibility for this note.

\section*{Acknowledgments}
A.D. acknowledges support from
the Croatian Science Foundation under the
project no. IP-2022-10-5008 (TEBAG),
the project ``Implementation of cutting-edge research and its application as part of the
Scientific Center of Excellence for Quantum and Complex Systems, and Representations of Lie Algebras'', Grant No. PK.1.1.10.0004, co-financed by the European
Union through the European Regional Development Fund -- Competitiveness and
Cohesion Programme 2021--2027,
and the European Union:
NextGenerationEU through the National Recovery and Resilience Plan 2021--2026, Institutional grant of University of Zagreb Faculty of Science (IK IA 1.1.3. Impact4Math).


\begin{thebibliography}{9}

\bibitem{Dtab} A.~Dujella, \emph{Infinite families of elliptic curves with high rank and prescribed torsion},
\url{https://web.math.pmf.unizg.hr/~duje/tors/generic.html}.

\bibitem{Er} Yu.~Eroshkin, personal communication to A.~Dujella (family of rank $4$ over $\Q(t)$ with torsion
$\Z/5\Z$, 2020).

\bibitem{J} M.~Jaswal, \emph{Elliptic curves with torsion group $\Z/6\Z$ and generic rank four}, preprint, Zenodo,
2026, \url{https://zenodo.org/records/23735797}.

\bibitem{HS} R.~Harron and A.~Snowden, \emph{Counting elliptic curves with prescribed torsion}, J.\ Reine Angew.\ Math.\
\textbf{729} (2017), 151--170.

\bibitem{PPVW} J.~Park, B.~Poonen, J.~Voight and M.~M.~Wood, \emph{A heuristic for boundedness of ranks of elliptic
curves}, J.\ Eur.\ Math.\ Soc.\ \textbf{21} (2019), 2859--2903.

\bibitem{MP} D.~Maulik and B.~Poonen, \emph{N\'eron--Severi groups under specialization}, Duke Math.\ J.\ \textbf{161}
(2012), 2167--2206.

\bibitem{SS} M.~Sch\"utt and T.~Shioda, \emph{Mordell--Weil Lattices}, Ergebnisse der Mathematik und ihrer Grenzgebiete
(3) \textbf{70}, Springer, 2019.

\bibitem{vL} R.~van Luijk, \emph{K3 surfaces with Picard number one and infinitely many rational points}, Algebra Number
Theory \textbf{1} (2007), 1--15.

\bibitem{Sil} J.~H.~Silverman, \emph{Heights and the specialization map for families of abelian varieties}, J.\ Reine
Angew.\ Math.\ \textbf{342} (1983), 197--211.

\end{thebibliography}
\end{document}